\documentclass{amsart}
\usepackage{amssymb}
\usepackage{lattice}
\usepackage{graphics}

\theoremstyle{plain}
\newtheorem{theorem}{Theorem}
\newtheorem{lemma}{Lemma}[section]
\newtheorem{proposition}[lemma]{Proposition}
\newtheorem{corollary}[lemma]{Corollary}
\newtheorem{example}[lemma]{Example}

\newtheorem*{stat}{\name}

\theoremstyle{definition}
\newtheorem{definition}[lemma]{Definition}

\newtheorem{remark}[lemma]{Remark}
\newtheorem{problem}{Problem}
\newcommand{\name}{testing}

\renewcommand{\gm}{\boldsymbol{\mu}}
\newcommand{\M}{\mathbf{M}}
\newcommand{\jz}{$\set{\jj, 0}$}

\newcommand{\dd}{\mathrm{d}}
\newcommand{\LL}{\E{L}}
\newcommand{\Nb}{\overline{\mathbb{N}}}

\newcommand{\incomp}{\parallel}
\newcommand{\FM}{\tup{F}_{\mathbf{M}}(J_1^4)}

\DeclareMathOperator{\J}{J}
\DeclareMathOperator{\Id}{Id}
\DeclareMathOperator{\Cl}{Cl}

\DeclareMathOperator{\Hom}{Hom}

\begin{document}

\title{The $M_{3}[D]$ construction and $n$-modularity}

\author{G.~Gr\"atzer}
\thanks{The research of the first author was supported by the
        NSERC of Canada.}
\address{Department of Mathematics\\
	  University of Manitoba\\
	  Winnipeg MN, R3T 2N2\\
	  Canada}
\email{gratzer@cc.umanitoba.ca}
\urladdr{http://www.maths.umanitoba.ca/homepages/gratzer.html/}

\author{F.~Wehrung}
\address{C.N.R.S.\\
   D\'epartement de Math\'ematiques\\
	  Universit\'e de Caen\\
	  14032 Caen Cedex\\
   France}
\email{wehrung@math.unicaen.fr}
\urladdr{http://www.math.unicaen.fr/\~{}wehrung}

\date{Sept. 17, 1998}

 \begin{abstract}
In 1968, E. T. Schmidt introduced the
$M_3[D]$ construction, an extension of the five-element
nondistributive lattice $M_3$ by a bounded distributive
lattice $D$, defined as the lattice of all triples $\vv<x, y,
z> \in D^3$ satisfying
$x \mm y = x \mm z = y \mm z$.  The lattice $M_3[D]$ is a
modular congruence-preserving extension of~$D$.

In this paper, we investigate this construction for an
arbitrary lattice
$L$.  For every $n > 0$, we exhibit an identity $\gm_n$ such
that $\gm_1$ is modularity and
$\gm_{n + 1}$ is properly weaker than $\gm_n$.  Let $\M_n$
denote the variety defined by $\gm_n$, the variety of
\emph{$n$-modular} lattices.  If $L$ is
$n$-modular, then $M_3[L]$ is a lattice, in fact, a
congruence-preserving extension of~$L$; we also prove that, in
this case, $\Id M_3[L] \iso M_3[\Id L]$.

We provide an example of a lattice $L$ such that $M_3[L]$ is
not a lattice.
This example also provides a negative
solution to a problem of R. W. Quacken\-bush: Is~the tensor
product $A \otimes B$ of two lattices $A$ and $B$ with zero
always a lattice.  We complement this result by generalizing the
$M_3[L]$ construction to an $M_4[L]$ construction. This yields, in
particular, a bounded modular lattice $L$ such that
$M_4 \otimes L$ is not a lattice, thus providing
a negative solution to Quacken\-bush's
problem in the variety $\M$ of modular lattices.

Finally, we sharpen a result of R. P. Dilworth: Every finite
distributive lattice can be represented as the congruence
lattice of a finite \emph{$3$-modular} lattice. We do this by
verifying that a construction of G.~Gr\"atzer, H. Lakser, and
E.~T. Schmidt yields a $3$-modular lattice.
 \end{abstract}

\maketitle

\section{Introduction}\label{S:Introduction} E. T. Schmidt
\cite{tS68} and \cite{tS74} introduced the following
construction. Let $M_3$ be the five-element, modular,
nondistributive lattice and let $D$ be a bounded distributive
lattice.  The lattice $M_3$ \emph{extended by} $D$, denoted by
$M_3[D]$, is the lattice of all triples $\vv<x, y, z> \in D^3$
satisfying
$x \mm y = x \mm z = y \mm z$; we call such triples
\emph{balanced}.  Then $M_3[D]$ is a (modular) lattice and
$M_3[D]$ and $D$ have isomorphic congruence lattices. Meet in
$M_3[D]$ is performed componentwise, while the join is the
smallest balanced triple in
$M_3[D]$ containing the triple formed by componentwise joins.

Note that the elements $\vv<x, 0, 0>$, $x \in D$, form a
sublattice of $M_3[D]$ isomorphic to $D$.  We identify $\vv<x,
0, 0> \in M_3[D]$ with $x \in D$, making
$M_3[D]$ an extension of $D$.

Let $L$ be a lattice.  A lattice $K$ is a
\emph{congruence-preserving extension} of
$L$, if $K$ is an extension of $L$ and every congruence of $L$
has
\emph{exactly one} extension to $K$.  Of course, then the
congruence lattice of $L$ is isomorphic to the congruence
lattice of~$K$. E.~T. Schmidt proved that \emph{$M_3[D]$ is a
congruence-preserving extension of $D$}.

This construction plays a central role in a number of papers
dealing with congruences of modular lattices, see G. Gr\"atzer
and E. T. Schmidt
\cite{GS95} and
\cite{GSa}, as two recent references.

This paper started with a problem proposed in G. Gr\"atzer and
E. T. Schmidt~\cite{GS95}: Does every lattice have a proper
congruence-preserving extension? (We solved this problem in
G. Gr\"atzer and F. Wehrung
\cite{GWa}.)  Of course, if the lattice is a bounded
distributive lattice $D$, then $M_3[D]$ is such an extension.  So
two problems were raised:

 1. For what classes of lattices $\mathbf{C}$, is $M_3[L]$ a
lattice for $L \in \mathbf{C}$?

 2. When is $M_3[L]$ a congruence-preserving extension of $L$?

Surprisingly, in addition to Schmidt's result ($M_3[D]$ is a
(modular) lattice provided that $D$ is a bounded distributive
lattice), we could only find one other relevant result in the
literature, see R.~W. Quackenbush \cite{rQ85}: if $L$ is
modular, then $M_3[L]$ is a lattice.

In Section~\ref{S:identities}, we define a lattice identity
$\gm_n$, for every $n > 0$, such that $\gm_1$ is equivalent to
the modular identity and $\gm_{n+1}$ is weaker than $\gm_n$. Let
$\M_n$ denote the variety defined by $\gm_n$; we call $\M_n$ the
variety of \emph{$n$-modular} lattices. We prove that if $L$ is
a $n$-modular lattice, then $M_3[L]$ is a lattice and it is a
congruence-preserving extension of $L$.
In~Section~\ref{S:m_n}, we verify that $\M_n \sci \M_{n + 1}$
and $\M_{n}$ gets very large as $n$ gets large:  $\UUm{\M_n}{n
< \go}$ generates the variety~$\mathbf{L}$.

We show, in Section~\ref{S:notlattice}, that $M_3[L]$ is, in
general, not a lattice.  In Section~\ref{S:nobounds}, we show
how we can remove the condition that $L$ be bounded in the
results of the previous sections. In Section~\ref{S:twoviews},
we explain how the two different definitions of $M_3[D]$ in the
literature can be reconciled using tensor products, and we
obtain the isomorphism $M_3[L] \iso M_3 \otimes L$, for any
lattice $L$ with zero which satisfies $\gm_n$ for some $n$.
It follows, then, that the result of
Section~\ref{S:notlattice} can be reinterpreted: there is a
lattice $L$ with zero such that $M_3 \otimes L$ is not a
lattice.  This solves, in the negative, a problem proposed in
R. W. Quackenbush \cite{rQ85}: Is the tensor product of two
lattices with zero always a lattice?  In fact, our
counterexample consists of two planar lattices.  In
Section~\ref{S:M4}, we show the there is a counterexample
consisting of two modular lattices, $M_4$ and the subspace
lattice of any infinite dimensional vector space. There is
another result on $n$-modular lattices in
Section~\ref{S:twoviews}: $\Id M_3[L] \iso M_3[\Id L]$.

In Section~\ref{S:Congruence}, we prove that every finite
distributive lattice can be represented as the congruence
lattice of a finite $3$-modular lattice $L$.  Without
$3$-modularity, this is a result of R. P. Dilworth.  We prove
this by verifying that the lattice~$L$ constructed by G.
Gr\"atzer, H. Lakser, and E. T. Schmidt \cite{AGLS95} to
represent $D$ is, in fact, $3$-modular.

The paper concludes with a discussion of some additional
results and a list of open problems in
Section~\ref{S:Discussion}.

 \section{The identities}\label{S:identities}

 Let $L$ be a lattice.  The triple $\vv<x, y, z> \in L^3$ is
\emph{balanced}, if
 \[
   x \mm y = x \mm z = y \mm z.
 \]
We denote by $M_3[L]$ the set of all balanced triples.  We
regard $M_3[L]$ as a subposet of $L^3$, in fact, a
meet-subsemilattice of $L^3$.

\begin{lemma}\label{L:closure}
 Let $L$ be a lattice. Then $M_3[L]$ is a lattice if{}f
$M_3[L]$ is a closure system in $L^3$.
 \end{lemma}

 \begin{proof}
 If $M_3[L]$ is a closure system and $\vv<x_0, y_0, z_0>$,
$\vv<x_1, y_1, z_1> \in M_3[L]$, then the closure of $\vv<x_0
\jj x_1, y_0 \jj y_1, z_0 \jj z_1>$ in $M_3[L]$ is the join of
$\vv<x_0, y_0, z_0>$ and $\vv<x_1, y_1, z_1>$.

Conversely, if $M_3[L]$ has joins, then the closure of
$\vv<x, y, z> \in M_3[L]$ is
 \[
   \vv<x, o, o> \jj \vv<o, y, o> \jj \vv<o, o, z>,
 \]
where $o$ is any element of $L$ contained in $x$, $y$, and $z$.

\vspace{-12pt}
 \end{proof}

Let us define the lattice polynomials $p_n$, $q_n$, and $r_n$,
for $n<\go$, in the variables $x$, $y$, and $z$:
 \begin{align*}
   p_0 &= x, &q_0 &= y, &r_0 &= z,\\
   p_{1} &= x \jj (y \mm z), &q_{1} &= y \jj (x
\mm z),
      &r_{1} &= z \jj (x \mm y),\\
   &&&\dots\\
   p_{n+1} &= p_n \jj (q_n \mm r_n), &q_{n+1} &= q_n \jj (p_n
\mm r_n),
      &r_{n+1} &= r_n \jj (p_n \mm q_n).
 \end{align*}

Let $\vv<x, y, z> \in L^3$.  Define, for $n > 0$,
 \[
   \vv<x, y, z>^{(n)} = \vv<p_n(x, y, z), q_n(x, y, z), r_n(x,
y, z)>
 \]
 Note that
 \begin{equation}\label{E:inreasing}
  \vv<x, y, z> \leq \vv<x, y, z>^{(1)} \leq \cdots \leq  \vv<x,
y, z>^{(n)}
\leq \cdots
 \end{equation}

\begin{definition}\label{D:identity}
 For $n > 0$, define the identity $\gm_n$ as $p_n = p_{n +
1}$.  Let $\M_n$ be the lattice variety defined by $\gm_n$. The
lattices in $\M_n$ are called \emph{$n$-modular}; lattices in
$\M_n - \M_{n-1}$ are called \emph{exactly $n$-modular} or
\emph{of modularity rank} $n$.  A lattice $L \nin \M_n$, for
all $n < \gw$, is of \emph{of modularity rank} $\infty$.
 \end{definition}

 \begin{lemma}\label{L:nstep}
 $L$ is an $n$-modular lattice if{}f, for all $a$, $b$, $c\in
L$, $\vv<a, b, c>^{(n)}$ is the closure of $\vv<a, b, c>$.
 \end{lemma}

 \begin{proof}
 This statement immediately follows from the definitions.
 \end{proof}

 \begin{corollary}\label{C:inclusions}
 The following inclusions hold:
 \[
   \M_1 \ci \M_2 \ci \cdots \ci \M_n \ci \cdots
 \]
 \end{corollary}

 \begin{corollary}\label{C:finite}
 For every finite lattice $L$, there is an integer $n > 0$ such
that $L$ is
$n$-modular.
 \end{corollary}

 \begin{proof}
 Indeed, if $L$ is finite, then the increasing sequence
$\vv<a, b, c>^{(n)}$ must terminate in $L^3$, so
Lemma~\ref{L:nstep} yields this result.
 \end{proof}

 \begin{corollary}\label{C:limit}
  $\UUm{\M_i}{i < \go}$ generates $\mathbf{L}$, the variety of all
lattices.
 \end{corollary}

 \begin{proof}
 Indeed, by Corollary~\ref{C:finite}, $\UUm{\M_i}{i < \go}$
contains all finite lattices and it is well-known that all
finite lattices generate the variety $\mathbf{L}$.
 \end{proof}

 \begin{lemma}\label{L:modular}
 $\M_1$ is the variety of modular lattices.
 \end{lemma}

 \begin{proof}
 If $L$ is a modular lattice, then computing in $\tup F_\M(3)$:
 \[
   p_2 = (x \jj (y \mm z)) \jj (y \jj (x \mm z)) \mm (z \jj (x \mm y)) =
       x \jj (y \mm z) = p_1.
 \]

Conversely, if $L$ is nonmodular, then it contains a pentagon $N_5 =
\set{o, a, b, c, i}$ (with zero $o$, unit $i$, and with $b < a$) as a
sublattice and
 \begin{align*}
  p_1(b, a, c) &=  b \jj (a \mm c) =  b,\\
  p_2(b, a, c) &=  (b \jj (a \mm c)) \jj ((a \jj (b \mm c))
 \mm (c \jj (a \mm b)) = a,
 \end{align*}
 so $\gm_1 \colon p_1 = p_2$ fails with $x = b$, $y = a$, and $z = c$.
 \end{proof}

On the other hand, $p_2 = p_3$ holds in $N_5$, so we obtain

 \begin{lemma}\label{L:N5}
 The variety $\mathbf{N}_5$ generated by $N_5$ is $2$-modular.
 \end{lemma}

E. T. Schmidt \cite{tS68} and \cite{tS74} proved that $M_3[L]$
is a modular lattice, if
$L$ is distributive.  We now prove the converse.

 \begin{lemma}\label{L:m3mod}
 $L$ be a lattice. Then $M_3[L]$ is a modular lattice
if{}f $L$ is distributive.
 \end{lemma}

 \begin{proof}
 So let $M_3[L]$ be modular.  Since $M_3[L]$ is an extension
of $L$, it follows that $L$ is modular.  If $L$ is not
distributive, then $L$ contains $M_3 = \set{o, a, b, c, i}$,
the five element modular nondistributive lattice, as a
sublattice.  Then
$\vv<o,o,o>$,
$\vv<a,a,a>$, $\vv<a,a,i>$, $\vv<b,c,o>$, $\vv<i,i,i> \in
M_3[L]$,
$\vv<a,a,i> \mm
\vv<b,c,o> = \vv<o,o,o>$, $\vv<a,a,a> \jj \vv<b,c,o> =
\vv<i,i,i>$, which easily imply that
 \[
   N_5 = \set{\vv<o,o,o>, \vv<a,a,a>, \vv<a,a,i>, \vv<b,c,o>,
\vv<i,i,i>}
 \] is the five-element nonmodular lattice, a sublattice of
$M_3[L]$, a contradiction.
 \end{proof}

The following lemma is due to E. T. Schmidt \cite{tS68}, for $n
= 1$, and to R.~W. Quackenbush \cite{rQ85}, for $n = 2$:

 \begin{lemma}\label{L:m3lattice}
   Let $n > 0$ and let $L \in \M_n$ be a lattice. Then
$M_3[L]$ is a lattice. Furthermore, if $L$ is bounded, then
$M_3[L]$ has a spanning $M_3$.
 \end{lemma}

 \begin{proof}
 By Lemma~\ref{L:closure}, we have to prove that $M_3[L]$ is a
closure system. For $\vv<x, y, z> \in M_3[L]$, define
 \begin{equation}\label{E:closure}
 \ol{\vv<x, y, z>} = \vv<p_n(x, y, z), q_n(x, y, z), r_n(x, y,
z)>.
 \end{equation}
 By \eqref{E:inreasing}, $\vv<x, y, z> \leq \ol{\vv<x, y,
z>}$.  Since a polynomial is isotone, $\vv<x, y, z> \leq
\vv<x', y', z'>$ implies that $\ol{\vv<x, y, z>} \leq
\ol{\vv<x', y', z'>}$.  Finally, let $\ol{\vv<x, y, z>} =
\vv<x^*, y^*, z^*>$; then
$p_n(x, y, z) = p_{n + 1}(x, y, z) = p_1(x^*, y^*, z^*)$, so
$\ol{\vv<x, y, z>}$ is closed.

The spanning $M_3$ is $\set{\vv<0, 0, 0>, \vv<1, 0, 0>, \vv<0,
1, 0>,
\vv<0, 0, 1>, \vv<1, 1, 1>}$.
 \end{proof}

 \begin{theorem}\label{T:congpres}
  Let $n > 0$ and let $L$ be a bounded $n$-modular lattice.
Then $M_3[L]$ is a lattice with a spanning $M_3$.  The map
 \[
   \ge \colon x \mapsto \vv<x, 0, 0>
 \]
 embeds $L$ into $M_3[L]$. If we identify $x \in L$ with $x\ge
= \vv<x, 0, 0> \in M_3[L]$, then the lattice $M_3[L]$ is a
congruence-preserving extension of $L$.
 \end{theorem}

\begin{proof}
 By Lemma~\ref{L:m3lattice}, $M_3[L]$ is a lattice.
Furthermore, $\ge$ is, obviously, an embedding; we identify
$x \in L$ with $x\ge = \vv<x, 0, 0>$.  Now let $\gQ$ be a
congruence of $L$. Form
$\gQ^3$, a congruence of $L^3$, and let $M_3[\gQ]$ be the
restriction of $\gQ^3$ to
$M_3[L]$.  We claim that $M_3[\gQ]$ is the unique extension of
the congruence $\gQ$ to
$M_3[L]$.  Since $M_3[\gQ]$ restricted to $L$ equals $\gQ$, it
is sufficient to prove the following two statements:
 \begin{enumerate}
 \item $M_3[\gQ]$ is a congruence of $M_3[L]$.
 \item  Every congruence $\gF$ of $M_3[L]$ is of the form
$M_3[\gQ]$, for some congruence
$\gQ$ of~$L$.
 \end{enumerate}

 \emph{Re:} (i). $M_3[\gQ]$ is obviously a meet-congruence on
$M_3[L]$.  It remains to prove the join substitution
property.  So let $\con \vv<x_0, y_0, z_0>=\vv<x_1, y_1,
z_1>(M_3[\gQ])$ and let $\vv<u, v, w> \in M_3[L]$.  Then $\con
x_0=x_1(\gQ)$ and so $\con x_0 \jj u=x_1 \jj u(\gQ)$.
Similarly, $\con y_0 \jj v=y_1 \jj v(\gQ)$ and
$\con z_0\jj w=z_1 \jj w(\gQ)$.  Since a polynomial has the
substitution property, we conclude that
 \[
   \con {p_n(x_0 \jj u, y_0 \jj v, z_0\jj w)}={p_n(x_1 \jj u,
y_1 \jj v,
      z_1 \jj w)}(\gQ),
 \]
 and similarly for $q_n$ and $r_n$. Thus
 \begin{align*}
   &\langle p_n(x_0 \jj u, y_0 \jj v, z_0 \jj w), q_n(x_0 \jj u,
      y_0 \jj v, z_0 \jj w),\\
 &\q\q r_n(x_0\jj u, y_0\jj v, z_0\jj w)\rangle\\
 & \equiv \langle p_n(x_1 \jj u, y_1 \jj v, z_1 \jj w), q_n(x_1
\jj u, y_1 \jj v, z_1 \jj w), \\
&\q\q r_n(x_0 \jj u, y_0 \jj v, z_0\jj w)\rangle
 \end{align*}
 modulo $\gQ^3$, and therefore, modulo $M_3[\gQ]$.
 By \eqref{E:closure}, this last congruence is the same as
 \[
   \con \vv<x_0, y_0, z_0> \jj \vv<u, v, w> =\vv<x_1, y_1,
     z_1> \jj \vv<u, v, w> (M_3[\gQ]),
 \]
 which was to be proved.

 \emph{Re:} (ii). Let $\gF$ be a congruence of $M_3[L]$ and
let $\gQ$ be the restriction of $\gF$ to~$L$.  We want to show
that $M_3[\gQ] = \gF$. Let $\vv<x_0, y_0, z_0>$,
$\vv<x_1, y_1, z_1> \in M_3[L]$.

 If $\con \vv<x_0, y_0, z_0>=\vv<x_1, y_1, z_1>(M_3[\gQ])$,
then $\con x_0=x_1(\gQ)$ in
$L$ and so  $\con \vv<x_0, 0, 0>=\vv<x_1, 0, 0>(M_3[\gQ])$.
Since
$M_3[\gQ]$ and $\gF$ agree on $L$, we conclude that
 \begin{equation}\label{E:first}
   \con \vv<x_0, 0, 0>=\vv<x_1, 0, 0>(\gF).
 \end{equation}
  Similarly, $\con \vv<y_0, 0, 0>=\vv<y_1, 0, 0>(\gF)$;
therefore,
\begin{align*}
 \vv<0, y_0, 0> &=   (\vv<y_0, 0, 0> \jj \vv<0, 0, 1>) \mm
\vv<0, 1, 0>\\
               &\equiv
   (\vv<y_1, 0, 0> \jj \vv<0, 0, 1>) \mm \vv<0, 1, 0> = \vv<0,
y_1, 0>
\pod{\gF}
 \end{align*}
 that is,
 \begin{equation}\label{E:second}
   \con \vv<0, y_0, 0>=\vv<0, y_1, 0>(\gF).
 \end{equation}
 Similarly,
 \begin{equation}\label{E:third}
   \con \vv<0, 0, z_0>=\vv<0, 0, z_1>(\gF).
 \end{equation} Joining the three congruences
\eqref{E:first}--\eqref{E:third}, we obtain that $\con
\vv<x_0, y_0, z_0>=\vv<x_1, y_1, z_1>(\gF)$.

Conversely, let $\con \vv<x_0, y_0, z_0>=\vv<x_1, y_1,
z_1>(\gF)$.  Meeting with $\vv<x_0
\jj x_1, 0, 0> \in M_3[L]$, we derive that $\con \vv<x_0, 0,
0>=\vv<x_1, 0, 0>(\gF)$ and so  $\con x_0=x_1(\gQ)$.
Similarly, $\con \vv<0, y_0, 0>=\vv<0, y_1, 0>(\gF)$.
Therefore,
 \[
    \vv<y_0, y_0, 1> = \con \vv<0, y_0, 0> \jj \vv<0, 0, 1>
=\vv<0, y_1,
      0> \jj \vv<0, 0, 1> = \vv<y_1, y_1, 1> (\gF)
 \]
 and meeting with $\vv<1, 0, 0> \in M_3[L]$, we conclude that
$\con
\vv<y_0, 0, 0>=\vv<y_1 ,0, 0>(\gF)$, that is, $\con
y_0=y_1(\gQ)$.  Similarly, $\con z_0=z_1(\gQ)$ and so $\con
\vv<x_0, y_0, z_0>=\vv<x_1, y_1, z_1>(\gQ^3)$, from which it
follows that
$\con \vv<x_0, y_0, z_0>=\vv<x_1, y_1, z_1>(M_3[\gQ])$.
 \end{proof}

\section{The variety $\M_n$}\label{S:m_n}
 In this section, we prove that $\M_{n}$ is properly contained
in $\M_{n + 1}$, for every
$n > 0$. It is obvious that $\M_{n} \ci \M_{n + 1}$.  To show
that the equality fails, we have to construct, for each $n >
0$, an exactly $(n+1)$-modular lattice $L_n$.  The lattice
$L_3$ is shown in Figure 1.  The definition of $L_n$ follows the
pattern of $L_3$, except that there are $n + 1$ $x$-s:
$x_0$, $x_1$, \dots, $x_n$; there are $n + 1$ $y$-s: $y_0$,
$y_1$, \dots,
$y_n$; and so there are $n + 1$ sublattices of the form $C_2^2$
in the middle.

\begin{figure}[bht]
\centerline{\includegraphics{Fig1.ill}}
\end{figure}

Since $L_n$ is planar and bounded, it is a lattice.

 \begin{theorem}\label{T:identity}
 $L_n$ is an exactly $(n + 1)$-modular lattice.
 \end{theorem}

 \begin{proof} It is easy to check that the free lattice on
$C_2 + C_1$ (see Figure VI.1.1 in \cite{GLT}) is $2$-modular.
This shows that $\gm_2$ holds in any lattice
at any triple $\vv<x,y,z>$ such that two of the
variables $x$, $y$ and $z$ are comparable.  So to check that
$\gm_2$ holds in $L_1$ at $\vv<x,y,z>$, we can assume that
$x$, $y$, and $z$ form an antichain; since there are very few
antichains of three elements in
$L_1$, it is very easy to compute that $\gm_2$ holds. So $L_1$
is (exactly) $2$-modular.

Now we induct on $n$. The interval $[x_1 \mm z_1, 1]$ of
$L_{n}$ is isomorphic to $L_{n - 1}$.  So if $x$, $y$, $z \in
[x_1 \mm z_1, 1]$, then $\gm_n$ holds at $\vv<x,y,z>$, therefore,
$\gm_{n + 1}$ holds at $\vv<x,y,z>$.  If two of $x$, $y$, $z$ are
in $[x_1 \mm z_1, 1]$, say, $y$, $z
\in [x_1 \mm z_1, 1]$, then replacing $x$ by
$\ol x = x \jj (y \mm z)$, and, similarly, for $y$ and $z$,
we have all three elements in
$[x_1 \mm z_1, 1]$ and $\gm_n$ holds for $\ol x$, $\ol y$,
$\ol z$, so $\gm_{n + 1}$ holds for
$x$, $y$,~$z$.  Since there is no three element antichain
outside of $[x_1
\mm z_1, 1]$, we are left with the case that two of $x$, $y$,
$z$ are not in $[x_1 \mm z_1, 1]$, say,
$x$ and $z$.  We cannot then have $x \leq x_1 \mm z_1$,
because there is no such antichain.  Similarly, $z \nleq x_1
\mm z_1$.  Theorefore, by symmetry, we can assume that $x = x_0$
and $z = z_0$.  It follows that $y \in [x_1 \mm z_1, y_0]$.  So
we have $p_1 = x_1$ and $q_1 = y$, $r_1 = z_1$, all in $[x_1 \mm
z_1, 1]$.  By induction,
$\gm_n$ holds for
$x_1$, $y$, and $z_1$ and so $\gm_{n + 1}$ holds for $x$,
$y$,~$z$.

It is clear that $\gm_n$ fails in $L_n$ with the substitution
$x = x_0$, $y = y_0$, and $z = z_0$ because $p_n(x_0, y_0,
z_0) = x_n < 1 = p_{n + 1}(x_0, y_0, z_0)$.
 \end{proof}

 \begin{corollary}\label{C:properinclusions}
 The following proper inclusions hold:
 \[
   \M_1 \sci \M_2 \sci \cdots \sci \M_n \sci \cdots
 \]
 \end{corollary}

\section{$M_3[L]$ is not always a lattice}\label{S:notlattice}
 In this section, we construct a bounded lattice $L$ such that
$M_3[L]$ is not a lattice.

 \begin{theorem}\label{T:nolattice}
 For the lattice $L$ of Figure 2, $M_3[L]$ is not a lattice.
Moreover, $L$ has modularity rank $\infty$.
 \end{theorem}

 \begin{proof}
 The reader can easily verify that $L$ is a lattice by
exhibiting the join- and meet-tables; for instance,
$x_i \mm z_j = c_{\min (i, j)}$, $x_i \jj z_j = 1$, and so on. By
Lemma~\ref{L:closure}, to show that $M_3[L]$ is not a lattice,
we have to verify that
$M_3[L]$ is not a closure system.  We claim that
$\vv<x_0, y_0, z_0>$ has no closure.  So let us assume to the
contrary that $\vv<\ol x, \ol y, \ol z>$ is the closure of
$\vv<x_0, y_0, z_0>$. Since $p_1(x_0, y_0, z_0) = x_1$,
$q_1(x_0, y_0, z_0) = y_0$,
$r_1(x_0, y_0, z_0) = z_1$, so by induction,
$\vv<\ol x, \ol y, \ol z>$ must contain all
$\vv<x_n, y_0, z_n>$, that is, $\vv<x_n, y_0, z_n> \leq
\vv<\ol x, \ol y, \ol z>$, for all
$n \geq 0$. On the other hand, $\vv<u_n, y_0, v_n>$ is
balanced, so
$\vv<\ol x, \ol y, \ol z> \leq \vv<u_n, y_0, v_n>$, for all
$n > 0$.  But there is no $\vv<\ol x, \ol y, \ol z> \in L$
satisfying $\vv<x_n, y_0, z_n> \leq
\vv<\ol x, \ol y, \ol z> \leq \vv<u_n, y_0, v_n>$, for all
$n > 0$, so $\vv<\ol x, \ol y, \ol z>$ does not exist.

\begin{figure}[thb]
 \centerline{\includegraphics{Fig2.ill}}
 \end{figure}

If $L$ was $n$-modular, for some $n < \go$, then $\vv<x_0, y_0,
z_0>^{(n)} = \vv<x_n, y_0, z_n>$ would be closed, but it is not.
 \end{proof}

We shall see in Section~\ref{S:twoviews} that $L$ provides a
negative solution to Quackenbush's problem, namely, $M_3 \otimes
L$ is not a lattice.

 \section{Removing the bounds}\label{S:nobounds}
 Most results of
Sections~\ref{S:identities}--\ref{S:notlattice} remain valid
without assuming that the lattice $L$ has a unit.  The only
exception is, of course, the statement that $M_3[L]$ has a
spanning~$M_3$.  If we do not assume that $L$ has a unit, then
the appropriate statement is that in $M_3[L]$, for every $a
\in M_3[L]$, there is a $i \in M_3[L]$ such that $(i]$ has a
spanning~$M_3$.

If we do not assume that $L$ has a zero, the definition of the
embedding
$\ge \colon x
\mapsto \vv<x, 0, 0>$ in Theorem~\ref{T:congpres} does not
make sense, affecting the crucial part about
con\-gru\-ence-preserving extensions. So we need to reformulate
Theorem~\ref{T:congpres}:

 \begin{theorem}\label{T:congpresgen}
   Let $n > 0$ and let $L$ be an $n$-modular lattice. Then
$M_3[L]$ is a lattice.  The map
 \[
   \gy \colon x \mapsto \vv<x, x, x>
 \]
 embeds $L$ into $M_3[L]$. If we identify $x \in L$ with
$x\gy = \vv<x, x, x> \in M_3[L]$, then the lattice $M_3[L]$ is a
congruence-preserving extension of $L$.
 \end{theorem}

 \begin{proof}
 The first part of the proof requires little change.

 Let $\gF$ be a congruence of $M_3[L]$ and let $\gQ$ be the
restriction of
$\gF$ to~$L$. We want to show that $M_3[\gQ] = \gF$. Let
$\vv<x_0, y_0, z_0>$, $\vv<x_1, y_1, z_1> \in M_3[L]$, and put
$o = \MMm{x_i \mm y_i}{i < 3}$.

 If $\con \vv<x_0, y_0, z_0> = \vv<x_1, y_1, z_1>(M_3[\gQ])$,
then $\con x_0=x_1(\gQ)$ in
$L$ and so  $\con \vv<x_0, x_0, x_0> = \vv<x_1, x_1, x_1>(\gF)$.
Therefore,
 \begin{align*}
   \vv<x_0, o, o> &= \vv<x_0, x_0, x_0> \mm \vv<x_0 \jj y_0,
o, o>\\
     &\equiv \vv<y_0, y_0, y_0> \mm \vv<x_0 \jj y_0, o, o> =
\vv<y_0, o, o>
     \pod{\gQ},
 \end{align*}
 that is,
 \begin{equation*}
    \vv<x_0, o, o> \equiv \vv<y_0, o, o> \pod{\gQ}.
 \end{equation*}

 Similarly, $\con \vv<o, y_0, o>=\vv<o, y_1, o>(\gF)$ and
$\con \vv<o, o, z_0>=\vv<o, o, z_1>(\gF)$.  Joining the three
congruences, we obtain $\con \vv<x_0, y_0, z_0>=\vv<x_1, y_1,
z_1>(\gF)$.

 The proof of the converse is similar to the original proof
with $o$ playing the role of
$0$ and $i = \JJm{x_i \jj y_i}{i < 3}$ playing the role of $1$.
 \end{proof}

 \section{Two views of $M_3[D]$}\label{S:twoviews}
 For a finite distributive lattice $D$, in the literature,
$M_3[D]$ is presented either as the lattice of balanced
triples $\vv<x, y, z> \in D^3$ (as we presented it in
Section~\ref{S:Introduction}) or as the lattice $M_3^P$, the
lattice of isotone maps from
$P = \tup J(D)$ (the poset of join-irreducible elements of
$D$) to $M_3$.  Either approach is convenient; both present a
modular lattice with a spanning $M_3$ with $D$ embedded as the
ideal generated by an atom of $M_3$ and the lattice is generated
by $M_3$ and $D$.  The second approach has the advantage that it
yields with no computation that
$M_3[D]$ is a modular lattice.  The~first approach, however,
better lends itself to generalization, as we did it in this
paper.

It follows from  A.~Mitchke and R. Wille \cite{MW73} that the
two constructions yield isomorphic lattices; indeed, both
constructions yield a modular lattice with a spanning
$M_3$ with $D$ embedded as the ideal generated by an atom of
$M_3$ and the lattice is generated by $M_3$ and $D$ and, up to
isomorphism, there is only one such lattice.

In this section we shall give a more direct explanation why
the two constructions yield isomorphic lattices.  To this end,
we introduce the concept of a capped tensor product from G.
Gr\"atzer and F. Wehrung \cite{GWb}.

\begin{definition}\label{D:tensor}
 Let $A$ and $B$ be \jz-semilattices. A \emph{bi-ideal} of
$A \times B$ is a subset $I$ of
$A \times B$ satisfying the following conditions:
 \begin{enumerate}
 \item $I$ is hereditary;
 \item $I$ contains $\nabla_{A,B} = (A \times \set{0})  \uu
(\set{0}
\times B)$;
 \item if $\vv<a_0, b>$, $\vv<a_1, b>\in I$, then $\vv<a_0 \jj
a_1, b> \in I$;
 \item if $\vv<a, b_0>$, $\vv<a, b_1>\in I$, then $\vv<a, b_0
\jj b_1> \in I$.
 \end{enumerate} \end{definition}

For $a \in A$ and $b \in B$, we define the bi-ideal
 \[
   a \otimes b =\nabla_{A,B}  \uu  \setm{\vv<x, y> \in A
\times B}{\vv<x, y> \leq
     \vv<a,b>}.
 \]

  The bi-ideal lattice of  $A \times B$ is an algebraic
lattice. The
\emph{tensor product}
$A \otimes B$ is the \jz-sub\-semi\-lat\-tice of compact
elements of the bi-ideal lattice of $A \times B$.

A bi-ideal $I$ is \emph{capped}, if there is a finite subset
$C$ of $A \times B$ such that
$I$ is the hereditary subset of $A \times B$ generated by $C$
along with $\nabla_{A, B}$. A tensor product $A \otimes B$ is
\emph{capped}, if all bi-ideals of $A \times B$ are
capped.  A capped tensor product is a lattice.

For a lattice $L$ with zero, let $L^-$ denote the
join-subsemilattice $L - \set{0}$.

Let $A \otimes B$ be a capped tensor product and let $I \in A
\otimes B$. We define a map
$\gf_I \colon A^-\to B$:

For $x \in A$, $x > 0$, let $\gf_I(x)$ be the largest element
$y$ in $B$ such that $\vv<x, y> \in I$.

 \begin{lemma}\label{L:map}
 $\gf_I$ maps $A^-$ into $B$ and
 \[
   \gf_I(x_0) \mm \gf_I(x_1) = \gf_I(x_0 \jj x_1),
 \]
 for $x_0$, $x_1 \in A^-$.
 \end{lemma}

 \begin{proof}
 First we show that $\gf_I(x)$ is defined, for all $x \in
A^-$.  Since $I$ is capped, we can write $I$ in the form $I =
\UUm{a_i \otimes b_i}{i < n} \uu \nabla_{A, B}$, where $n$ is a
natural number, $a_i \in A$, $b_i \in B$, for $i < n$.  Now
define
 \[
   y_i =
       \begin{cases}
          b_i, &\text{if $x \leq a_i$;}\\
          0,   &\text{otherwise,}
       \end{cases}
 \]
 for $i < n$ and let $y = \JJm{y_i}{i < n}$.  By definition,
$\vv<x, y_i>
\in I$, so by
\ref{D:tensor}(iv), $\vv<x, y> \in I$.  Now let $\vv<x, z> \in
I$, for some
$z \in B$. Then $\vv<x, z> \in a_i \otimes b_i$, for some $i <
n$, and so $z \leq b_i
\leq y$.  This proves that $y$ satisfies the requirements in
the definition of~$\gf_I(x)$.

Now $\gf_I(x_0 \jj x_1) \leq \gf_I(x_0)$ is obvious, hence,
$\gf_I(x_0 \jj x_1) \leq
\gf_I(x_0) \mm \gf_I(x_1)$.  Conversely, $\gf_I(x_0) \mm
\gf_I(x_1) \leq
\gf_I(x_0)$, so
$\vv<x_0, \gf_I(x_0) \mm \gf_I(x_1)> \in I$; similarly,
$\vv<x_1,
\gf_I(x_0) \mm
\gf_I(x_1)> \in I$, therefore, by \ref{D:tensor}(iii),
$\vv<x_0 \jj x_1,
\gf_I(x_0) \mm
\gf_I(x_1)> \in I$, implying that $\gf_I(x_0) \mm \gf_I(x_1)
\leq \gf_I(x_0
\jj x_1)$.
 \end{proof}

Let $B^\dd$ denote the dual lattice of $B$ and let
$\Hom_{\jj}(A^-, B^\dd)$ denote the lattice of
join-homomorphisms from $A^-$ to $B^\dd$, ordered
componentwise, as a subset of
$B^{A^-}$ (not $(B^\dd)^{A^-}$).

 \begin{theorem}\label{T:repr}
 Let $A$ and $B$ be lattices with zero and let $A \otimes B$ be
capped. Then $\ge \colon I \mapsto \gf_I$ defines an
isomorphism between $A \otimes B$ and $\Hom_{\jj}(A^-,B^\dd)$.
 \end{theorem}

 \begin{proof}
 Lemma~\ref{L:map} states that the map is well-defined.  Since
$\vv<x, y>
\in I$ if{}f $y
\leq \gf_I(x)$, it follows that $\gf_I$ determines $I$ and so
$\ge$ is one-to-one.  To show that $\ge$ is onto, let $\gf \in
\Hom_{\jj} (A^-, B^\dd)$ and define $I =
\setm{\vv<x, y> \in A \times B}{y \leq \gf(x)}$.  Since $\gf$
is a join-homomorphism, it follows that $I$ is a bi-ideal and
$\gf = \gf_I$.
 \end{proof}

\begin{corollary}\label{C:M3otLlatt}
 Let $L$ be a lattice with zero. Then the following
conditions are equivalent:
 \begin{enumerate}
 \item $M_3\otimes L$ is a lattice.
 \item For all $\vv<x, y, z>\in L^3$, there exists $n>0$ such
that $\vv<x, y, z>^{(n)}=\vv<x, y, z>^{(n + 1)}$.
 \end{enumerate}
Furthermore, if \tup{(i)} is satisfied, then
$M_3\otimes L\iso M_3[L]$.

In particular, if $L$ is $n$-modular, for some $n$, then $M_3
\otimes L$ is a lattice, and
$M_3 \otimes L \iso M_3[L]$.
 \end{corollary}

\begin{proof}
 Since $M_3$ is finite, $M_3\otimes L$ is a lattice if{}f
$M_3\otimes L$ is a capped tensor product, see Theorem~3 of
\cite{GWc}. Furthermore, in the same theorem, it is stated that
this is equivalent to saying that, for every antitone map
$\gx\colon\J(M_3)\to L$, the adjustment sequence of $\gx$
terminates after a finite number of steps. Here,
$\J(M_3)=\set{a,b,c}$, and the ordering on $\J(M_3)$ is
trivial, thus every map from
$\J(M_3)$ to $L$ is antitone. Identify $\gx\colon\J(M_3)\to L$
with the triple
$\vv<\gx(a),\gx(b),\gx(c)>$. With this identification, the
adjustment sequence of $\gx$ is easily seen to be the sequence
of all $\vv<\gx(a),\gx(b),\gx(c)>^{(n)}$, $n > 0$. The
equivalence between (i) and (ii) follows.

Now assume that $M_3\otimes L$ is a lattice. Again,
$M_3\otimes L$ is a capped tensor product, thus, by
Theorem~\ref{T:repr}, $M_3\otimes
L\iso\Hom_{\jj}(M_3^-,L^\dd)$. For all
$\gx\in\Hom_{\jj}(M_3^-,L^\dd)$, we can identify $\gx$ with
the triple
$\vv<\gx(a),\gx(b),\gx(c)>$, which, by Lemma~\ref{L:map}, is a
balanced triple. The isomorphism $M_3 \otimes L \iso M_3[L]$
follows.
 \end{proof}

Now the solution of Quackenbush's problem (discussed in the
Introduction) easily follows:

 \begin{corollary}\label{C:Quackenbush}
 Let $L$ be the lattice of Theorem~\ref{T:nolattice}.  Then $M_3
\otimes L$ is not a lattice.
 \end{corollary}
 \begin{proof}
 This is obvious by Theorem~\ref{T:nolattice} and
Corollary~\ref{C:M3otLlatt}.\
 \end{proof}

 \begin{corollary}\label{C:twofaces}
  Let $L$ be a lattice with zero.  If $L$ is $n$-modular, for
some $n$, then
 \begin{equation}
   M_3 \otimes L \iso M_3[L].\tag{i}
 \end{equation}
 If $L$ is a finite distributive lattice and $P = \tup J(L)$,
the poset of join-irreducible elements of $D$, then
 \begin{equation}
   M_3 \otimes L \iso M_3^P.\tag{ii}
 \end{equation}
 \end{corollary}

 \begin{proof} Part (i) follows immediately from
Corollary~\ref{C:M3otLlatt}.

If $L$ is a finite distributive lattice, then, by
Theorem~\ref{T:repr},
$M_3 \otimes L$ is isomorphic to $\Hom_{\jj}(L^-,M_3)$, and
any $\gf\in\Hom_{\jj}(L^-,M_3)$ can be identified with an
isotone map from $P$ into $M_3$.
 \end{proof}

Combining (i) and (ii) of Corollary~\ref{C:twofaces}, we
obtain the desired isomorphism:

 \begin{corollary}\label{C:finaliso}
  Let $D$ be a finite distributive lattice.  Then
 \[
   M_3^P \iso M_3[D],
 \]
 where $P = \tup J(D)$, the poset of join-irreducible elements
of $D$.
 \end{corollary}

For a given lattice, $n$-modularity has the following
algebraic meaning:

 \begin{proposition}\label{P:M3Lgo}
 Let $L$ be a lattice with zero. Then the following
conditions are equivalent:
 \begin{enumerate}
 \item $M_3\otimes L^\go$ is a lattice.
 \item $L$ is $n$-modular, for some $n > 0$.
 \end{enumerate}
 \end{proposition}

\begin{proof}\hfill

 (ii) \emph{implies} (i). If $L$ is $n$-modular, then $L^\go$ is
$n$-modular, thus, by Corollary~\ref{C:M3otLlatt}, $M_3\otimes
L^\go$ is a lattice.

(i) \emph{implies} (ii). Let us assume that the modularity rank
of $L$ is $\infty$. For all $n > 0$, there exists, by
definition, a triple $\vv<x_n, y_n, z_n>\in L^3$ such that
 \[
   \vv<x_n, y_n, z_n>^{(n)}<\vv<x_n, y_n, z_n>^{(n + 1)}.
 \]
 Define the following elements of $L^\go$:
 \begin{align*}
     x &= \vv<{x_n \mid n > 0}>,\\
     y &=\vv<{y_n \mid n > 0}>,\\
     z &= \vv<{z_n \mid n > 0}>.
 \end{align*}
 Then, $\vv<x,y,z>^{(n)}<\vv<x,y,z>^{(n+1)}$ holds, for all
$n > 0$. By Corollary~\ref{C:M3otLlatt}, $M_3\otimes L^\go$ is
not a lattice.
 \end{proof}

Another algebraic consequence of $n$-modularity is the
following:

 \begin{proposition}\label{P:ideals}
 Let $L$ be an $n$-modular lattice, for some $n > 0$.  Then
 \[
   M_3[\Id L] \iso \Id M_3[L].
 \]
 \end{proposition}

 \begin{proof}
 Let $U \in M_3[\Id L]$, that is, $U = \vv<I, J, K> \in M_3[\Id
L]$.  We define
 \[
  \widehat U = \setm{\vv<a, b, c> \in M_3[L]}{a \in I,\ b \in
J,\ c \in K},
 \]
 that is, $\widehat U = (I \times J \times K) \ii M_3[L]$.
Then $\widehat U$ is obviously a hereditary subset of
$M_3[L]$.  We claim that it is join closed.  Indeed, let
$\vv<a_0, b_0, c_0>$, $\vv<a_1, b_1, c_1> \in \widehat U$, set
 \[
   \vv<a_0 \jj a_1, b_0 \jj b_1, c_0 \jj c_1> = \vv<x, y, z> \in
     L^3,
 \]
 and define
$\vv<x, y, z>^{(n)}$ as in Section~\ref{S:identities}. We have
$x \in I$, $y \in J$,
$z \in K$, so $y \mm z \in J \mm K \ci I$ and
$x \jj (y \mm z) \in I$; similarly,
$y \jj (x \mm z) \in J$, $z \jj (x \mm y) \in K$.
Thus
 \[
   \vv<x, y, z>^{(1)} = \vv<x \jj (y \mm z), y \jj (x \mm z), z
      \jj (x \mm y)> \in I \times J \times K.
 \]
 By~induction, $\vv<x, y, z>^{(n)}\in I \times J \times K$.
Since $L$ satisfies $\gm_n$, $\vv<a_0, b_0, c_0> \jj \vv<a_1,
b_1, c_1> = \vv<x, y, z>^{(n)}$, so $\vv<a_0, b_0, c_0> \jj
\vv<a_1, b_1, c_1> \in I \times J \times K$, and therefore,
$\vv<a_0, b_0, c_0> \jj \vv<a_1, b_1, c_1> \in \widehat U$,
proving that $\widehat U$ is an ideal of $M_3[L]$.

Since every $a \in I$ can be augmented to an $\vv<a, b, c> \in
M_3[L]$ with suitable $b
\in J$, $c \in K$ (and similarly for  $b \in J$ and for $c \in
K$), it follows that
$\widehat U$ determines $U$.  Therefore, $\gf \colon U \mapsto
\widehat U$ is a one-to-one map from $M_3[\Id L]$ into $\Id
M_3[L]$.

To complete the proof, we have to prove that $\gf$ is onto. So
let $X \in
\Id M_3[L]$. Define
\[
   I = \setm{a \in L}{\vv<a, b, c> \in X, \text{ for some } b,\
     c \in M_3[L]}.
 \]
 Since $X$ is hereditary,
so is $I$.  Now let $a_0$, $a_1 \in I$.  Then
$\vv<a_0, b_0, c_0>
\in X$ and $\vv<a_1, b_1, c_1> \in X$, for some $b_0$, $b_1$,
$c_0$, $c_1
\in M_3[L]$. Let
\[
   \vv<a_0, b_0, c_0> \jj \vv<a_1, b_1, c_1> = \vv<x, y, z> \in
     M_3[L].
 \]
 Then $x \in I$ and
$a_0 \jj a_1 \leq x$, hence, $a_0 \jj a_1 \in I$, proving that
$I$ is an ideal. Similarly, one can define the ideals $J$ and
$K$ of $L$, by permuting the three coordinates in $L^3$.

Now we prove that
 \[
   \vv<I, J, K> \in M_3[\Id L].
 \]
 Indeed, let $w \in J \mm K$. Then $w \in J$, so there exist $i
\in I$ and $k \in K$ such that $\vv<i, w, k> \in X$. Similarly,
since $w\in K$, there exist $i'\in I$ and $j\in J$ such that
$\vv<i',j,w>\in X$. Put
$o=i\mm i'\mm j\mm k\mm w$. Since $X$ is an ideal, we have
$\vv<o, w, o>
\in X$. Similarly, $\vv<o, o, w> \in X$.  Since $X$ is join
closed,
 \[
   \vv<o, w, o> \jj \vv<o, o, w> = \ol{\vv<o, w, w>} = \vv<w, w,
      w> \in X,
 \]
 therefore, we conclude that $w \in I$, so $J \mm K
\ci I$. Similarly, $I \mm K \ci J$ and $I \mm J \ci K$, so
 \[
   I \mm K = I \mm J = J \mm K,
 \]
 that is, $\vv<I, J, K> \in M_3[\Id L]$.

Finally, we prove that, for $U = \vv<I, J, K>$, we have
$\widehat U = X$. Indeed, by the definitions of $I$, $J$, and
$K$, it is obvious that $X \ci \widehat U$. So let $\vv<x, y,
z> \in \widehat U$, that is, $\vv<x, y, z> \in M_3[L]$ and $x
\in I$, $y
\in J$, $z \in K$.  Since $x \in I$, there exist $j \in J$ and
$k \in K$ such that
$\vv<x,j,k>\in X$. Similarly, there are $i\in I$ and $k'\in K$
such that $\vv<i,y,k'>\in X$, and there are
$i'\in I$ and $j'\in J$ such that $\vv<i',j',z>\in X$.
Therefore,
 \[
   \vv<x, y, z> \leq \vv<x, j, k> \jj \vv<i, y, k'> \jj \vv<i',
     j', z> \in X,
 \]
proving that $\widehat U \ci X$.
 \end{proof}

\begin{remark}
E.~T. Schmidt suggested that we consider, for an arbitrary lattice
$L$ with zero, the least join-congruence $\gQ$ of $L^3$
identifying all triples $\vv<0, x, x>$, $\vv<x, 0, x>$, and
$\vv<x, x, 0>$, for all $x \in L$, and relate the quotient lattice
$L/\gQ$ to $M_3[L]$.

An alternative description of $\gQ$ is the following. For all
$\vv<a,b,c>\in L^3$, denote by $\Cl\vv<a,b,c>$ the least ideal
$I$ of $L^3$ containing $\vv<a,b,c>$ and such that if
$\vv<x,y,z>\in I$, then $\vv<x,y,z>^{(1)}\in I$. In particular,
if $L$ is $n$-modular, for some $n$, then $\Cl\vv<x,y,z>$ is just
the ideal of $L^3$ generated by the closure $\ol{\vv<x,y,z>}$ of
$\vv<x,y,z>$. Then it is not hard to verify that
 \[
 \con \vv<x_0,y_0,z_0>=\vv<x_1,y_1,z_1>(\gQ)
 \qq\text{if{f}}\qq
 \Cl\vv<x_0,y_0,z_0>=\Cl\vv<x_1,y_1,z_1>,
 \]
for all $x_i$, $y_i$, $z_i$ in $L$ and for $i<2$. Hence, if $L$
is $n$-modular, for some $n$, then $L^3/\gQ$ is isomorphic to
$M_3[L]$. Furthermore, $\Cl\vv<x_0,y_0,z_0>=\Cl\vv<x_1,y_1,z_1>$
is also equivalent to
 \[
 (a\otimes x_0)\jj(b\otimes y_0)\jj(c\otimes z_0)=
 (a\otimes x_1)\jj(b\otimes y_1)\jj(c\otimes z_1),
 \]
thus $L^3/\gQ$ is isomorphic to $M_3\otimes L$, in general.
\end{remark}

\section{The $M_4[L]$ construction}\label{S:M4}

Let $M_4$ denote the lattice of height $2$ with four atoms,
$a$, $b$, $c$, and $d$.

In this section, we prove that $M_4 \otimes L$ is not a
lattice, for a suitable modular lattice $L$ with zero, thereby
showing that R. W. Quackenbush's problem (discussed in the
Introduction) has a negative solution also for modular
lattices.  We also find new examples of nonmodular tensor
products that are not lattices, for instance, $(M_3 \otimes M_3)
\otimes L$.

The $M_3[L]$ construction has a natural extension to $M_4$.
For every lattice $L$, define
 \[
   M_4[L]=\setm{\vv<x_0 ,x_1, x_2, x_3> \in L^4}{x_i \mm x_j =
     x_0\mm x_1, \text{ for } i \ne j}.
 \]
 As in Lemma~\ref{L:closure}, it is easy to prove that
$M_4[L]$ is a meet-subsemilattice of $L^4$, and that it is a
lattice if and only if it is a closure system in $L^4$.

Just as for triples, define the lattice polynomials $q_0$,
$q_1$, $q_2$, and $q_3$ in four variables, $x_0$, $x_1$, $x_2$,
and $x_3$, as follows:
 \begin{gather*}
   q_{i0}=x_i,\q \text{ for } i < 4,\\
   q_{0, n + 1} = q_{0n} \jj (q_{1n} \mm q_{2n}) \jj (q_{1n}
\mm q_{3n})
\jj (q_{2n} \mm
     q_{3n}),
 \end{gather*} and, cyclically, define $q_{1, n + 1}$, $q_{2,
n + 1}$, and $q_{3, n + 1}$.

Let $\vv<x_0, x_1, x_2, x_3>\in L^4$. Define, for $n<\go$,
 \begin{align}\label{Eq:Adj4}
   &\vv<x_0, x_1, x_2, x_3>^{(n)}\\
   &=\vv<q_{0n}(x_0, x_1, x_2, x_3),\dots,q_{3n}(x_0, x_1, x_2,
     x_3)>.\notag
 \end{align}

The proof of the following result is very similar to the proof
of Corollary~\ref{C:M3otLlatt}, thus we will omit it.

\begin{proposition}\label{P:M4otLlatt}
 Let $L$ be a lattice with zero. Then the following are
equivalent:
 \begin{enumerate}
 \item $M_4\otimes L$ is a lattice.
 \item For all $\vv<x_0, x_1, x_2, x_3> \in L^4$, there exists
$n > 0$ such that
 \[
  \vv<x_0, x_1, x_2, x_3>^{(n)} = \vv<x_0, x_1, x_2, x_3>^{(n +
      1)}.
 \]
 \end{enumerate}

Furthermore, if \tup{(i)} is satisfied, then $M_4 \otimes L
\iso M_4[L]$.
\end{proposition}

Again, it is not difficult to verify that if $L$ is
distributive, then it satisfies the identity $\vv<x_0, x_1,
x_2, x_3>^{(2)} = \vv<x_0, x_1, x_2, x_3>^{(1)}$. However,
this is no longer true for modular lattices, as witnessed by
the main result of this section:

\begin{theorem}\label{T:DHWexple}
 Let $V$ be an infinite dimensional vector space over a field
$K$. Denote by $\LL(V)$ the lattice of all subspaces of
$V$. Then $M_4 \otimes \LL(V)$ is not a lattice.
 \end{theorem}

Note the contrast with the $M_3$ case: $\LL(V)$ is a
\emph{modular} lattice, so
$M_3 \otimes \LL(V)$ is a lattice.

 \begin{proof} We shall work with the lattice $\FM$ introduced
in A. Day, C. Herrmann, and R.~Wille \cite{DHW} as follows. By
definition, $\FM$ is the modular lattice generated by the
elements $a$, $b$, $c$, and $d$, subject to the relations
 \begin{gather*}
    a \mm b = a \mm c = a \mm d = b \mm c = b \mm d = c \mm d =
0,\\
    a \jj b = a \jj c = a \jj d = b \jj d = c \jj d = 1.
 \end{gather*}
 It is proved in \cite{DHW} that $\FM$ is isomorphic to a
sublattice of the lattice of all subgroups of a free abelian
group of infinite rank. Replacing the free Abelian group on a
countably infinite number of generators by a vector space $U$
of countably infinite dimension over a field $K$, it is easy
to see that the same construction shows that $\FM$ is
isomorphic to a sublattice of the subspace lattice $\LL(U)$.

The classical dualization map
 \[
   X \mapsto X^{\bot} = \setm{\gx \in U^*}{\gx[X] = \set{0}},
 \]
 where $U^*$ denotes the dual of $U$, is a dual embedding from
$\LL(U)$ into $\LL(U^*)$. Therefore, the dual lattice
$L = \FM^\dd$ of
$\FM$ embeds into $\LL(U^*)$. Now let $V_0$ be a vector
subspace of $U^*$ satisfying the following properties:

 \begin{enumerate}
 \item $V_0$ has countably infinite dimension.
 \item $(X + Y) \ii V_0 = (X \ii V_0) + (Y \ii V_0)$, for all
$X$, $Y\in L$.
 \item $(X - Y) \ii V_0 \ne \es$, for all $X$, $Y$ in $L$ such
that $X \nci Y$.
 \end{enumerate} Then $X\mapsto X \ii V_0$ is a lattice
embedding from $L$ into $\LL(V_0)$. Hence we have reached
the following conclusion:

 \medspace

\noindent \emph{$L$ embeds into $\LL(V_0)$, where $V_0$ is a
countably infinite dimensional vector space over~$K$}.

 \medspace

$\LL(V_0)$ embeds into $\LL(V)$, thus $L$ embeds into $\LL(V)$.

Now suppose that $M_4 \otimes \LL(V)$ is a lattice. By
Proposition~\ref{P:M4otLlatt}, the adjustment sequence
(\ref{Eq:Adj4}) based on any quadruple of elements of $\LL(V)$
terminates. Thus, \emph{a~fortiori}, the same holds for
quadruples of elements of $L$.

However, we shall now prove that there exists a quadruple of
elements of
$L$ whose adjustment sequence does not terminate, thus
completing the proof. For this, we need the following
description of $L$, obtained by dualizing the one given in
\cite{DHW} for $\FM$.

Put $\Nb = \go \uu \set{\infty}$. For $\vv<i, j>$ and $\vv<k,
l>$ in $\Nb
\times \Nb$, let us write
$\vv<i, j> \sim \vv<k, l>$, if $x \equiv y\pmod{2}$, for all
$x$, $y \in
\set{i, j, k, l} - \{\infty\}$. Then, we have
 \[
   L = \setm{\vv<i, j> \in \Nb\times\Nb}{\vv<i,
j>\sim\vv<\infty, \infty>}.
 \]
 The least element of $L$ is $\vv<0, 0>$, the largest element
of $L$ is
$\vv<\infty, \infty>$.

Denote by $\mm$ and $\jj$ the infimum and the supremum on $\Nb$.
The meet and the join of $L$ are given as follows:
 \begin{align}
 \vv<i, j> \mm \vv<k, l> =
				\begin{cases}
      \vv<i \mm  k, j \mm l>,\\
     \qq \text{ if } \vv<i, j> \sim \vv<k, l>;\\
      \vv<(i - 1) \mm (j - 1)\mm k, (i - 1) \mm (j - 1) \mm
l>,\\
      \qq \text{ if } \vv<i, j> \not \sim \vv<k, l> \text{ and
}i \mm j\geq k \mm  l;\\
      \vv<i \mm (k - 1) \mm (l - 1), j \mm (k - 1) \mm (l -
1)>,\\
      \qq \text{ if } \vv<i, j> \not\sim \vv<k, l> \text{ and
}i \mm  j<k
\mm  l;
				\end{cases}\label{Eq:meet}\\
  \vv<i, j> \jj \vv<k, l> =
				\begin{cases}
       \vv<i\jj k, j\jj l>,\\
       \qq\text{ if } \vv<i, j> \sim \vv<k, l>;\\
       \vv<(i + 1)\jj(j + 1)\jj k, (i + 1)\jj(j + 1)\jj l>,\\
       \qq\text{ if } \vv<i, j>\not \sim \vv<k, l> \text{ and
}i\jj j\leq k\jj l;\\
       \vv<i\jj(k + 1)\jj(l + 1), j\jj(k + 1)\jj(l + 1)>,\\
       \qq\text{ if } \vv<i, j>\not \sim \vv<k, l> \text{ and
}i\jj j>k\jj l.
				\end{cases}
 \label{Eq:join}
 \end{align} The base quadruple $ \vv<x, y, z, t>$ of elements
of $L$ is given by
 \[
    x = \vv<0, \infty>;\q y = \vv<1, \infty>;\q
    z = \vv<\infty, 1>;\q t = \vv<\infty, 0>.
 \] For all $n > 0$, put
 \[
    \vv<x, y, z, t>^{(n)} = \vv<x^{(n)}, y^{(n)}, z^{(n)},
      t^{(n)}>.
 \]
 Then an easy (though somehow tedious) induction proof,
based on the formulas (\ref{Eq:meet}) and (\ref{Eq:join}),
gives that for all $n > 0$, we have
 \begin{alignat*}{2}
    x^{(2n + 1)} &= x^{(2n + 2)} &&= \vv<2n + 2, \infty>,\\
    t^{(2n + 1)} &= t^{(2n + 2)} &&= \vv<\infty, 2n + 2>,\\
        y^{(2n)} &= y^{(2n + 1)} &&= \vv<2n + 1, \infty>,\\
         z^{(2n)}&= z^{(2n + 1)} &&= \vv<\infty, 2n + 1>.
 \end{alignat*}
 In particular, the sequence $\vv<x, y, z, t>^{(n)}$, $n > 0$,
is not eventually constant.
 \end{proof}

\begin{corollary}\label{C:M3ofM3}
 Let $V$ be an infinite dimensional vector space over a field
$K$.  Then $M_3 \otimes M_3[\LL(V)]$ is not a lattice.
 \end{corollary}

 \begin{proof}
 Define $K = M_3[\LL(V)]$ and $L = \LL(V)$. Since $L$ is
modular, it follows from Corollary~\ref{C:M3otLlatt} that
$M_3[L]\iso M_3\otimes L$. By Proposition~2.9 of \cite{GWb},
the tensor product of semilattices with zero is associative,
thus we have
 \[
   M_3\otimes(M_3\otimes L)\iso(M_3\otimes M_3)\otimes L.
 \]
 Thus, in order to prove that $M_3\otimes K$ is not a lattice,
it suffices to prove that
$(M_3\otimes M_3)\otimes L$ is not a lattice.

 Now, we note that the following four elements
 \[
    t = \vv<1, 0, 0>,\quad u = \vv<0, a, b>,\quad
    v = \vv<0, b, c>,\quad w = \vv<0, c, a>
 \]
 of $M_3[M_3]$ have pairwise meet $\vv<0, 0, 0>$ and pairwise
join $\vv<1, 1, 1>$, thus they generate a $0$-sublattice
isomorphic to $M_4$. Hence, there exists a zero preserving
embedding $f$ of $M_4$ into $M_3\otimes M_3 \iso M_3[M_3]$.

Now we need a very special case of Corollary 3.8 of \cite{GWb}:

\smallskip

 \emph{Let $A$, $A'$, $B$  be lattices with zero such that $A$
is a
$\set{0}$-sublattice of
$A'$.  If $A' \otimes B$ is a lattice, then $A \otimes B$ is a
lattice.}

\smallskip

Apply this with $A = M_4$, $A' = M_3[M_3]$, and $B = L$.  By
Theorem~\ref{T:DHWexple}, $M_4\otimes L$ is not a lattice.
Therefore, by the above statement, $(M_3\otimes M_3)\otimes L$
is not a lattice either.
 \end{proof}

\begin{corollary}\label{C:nonmod} There is a modular lattice
$M$ such that $M_3[M]$ is of modularity rank $\infty$.
 \end{corollary}

\begin{proof}
 We can take $M = \LL(V)$, for any infinite dimensional vector
space $V$.  Indeed, if $M_3[M]$ is $n$-modular, for some $n > 0$, then,
by Theorem~\ref{T:congpres}, $M_3 \otimes M_3[M]$ would be a lattice,
contradicting Corollary~\ref{C:M3ofM3}.
 \end{proof}

\section{Congruence lattices}\label{S:Congruence}
 In this section, we prove that every finite distributive lattice can be
represented as the congruence lattice of a finite $3$-modular lattice
$L$.

 In  G. Gr\"atzer, H. Lakser, and E. T. Schmidt \cite{AGLS95},
it is proved that every finite distributive lattice $D$ can be
represented as the congruence lattice of a finite planar
lattice $L$. This lattice $L$ has the following properties:
 \begin{enumerate}
 \item[(C1)] $L$ has a $\set{0, 1}$-sublattice, $G$, the
\emph{grid}, which is of the form $C \times D$, where $C$ and
$D$ are finite chains.
 \item[(C2)]  Every element of $H  =  L - G$ is doubly
irreducible in $L$.
  \end{enumerate}

It follows from (C1) that for every element $x \in L$, there is
a largest grid element $\ul x$ with $\ul x \leq x$; dually,
there is a smallest grid element $\ol x \in G$ with $x \leq \ol
x$.
 \begin{enumerate}
 \item[(C3)] For $x$, $y \in H$, if $\ul x = \ul y$, then $x =
y$; and dually.
 \item[(C4)] For every $x \in H$, either the interval $[\ul x,
\ol x]_G$ is a prime interval and $[\ul x, \ol x]_L$ is the
three-element chain, or $[\ul x, \ol x]_G$ is a prime square
and $[\ul x, \ol x]_L$ is an $M_3$.
  \end{enumerate} Recall that a \emph{prime square} is an
interval of length two isomorphic to $C_2^2$.  For a grid
element $x$, we shall use the notation $\vv<x_C, x_D>$, where
$x_C \in C$ and $x_D \in D$.

By (C4), if $[a, b]$ is a prime interval in $G$, then $[a, b]_L
- [a, b]_G$ is either empty or it is a singleton; in the latter
case, we denote the new element by $n(a, b)$. Similarly, if
$[a, b]$ is a prime square in $G$, then the set $[a, b]_L -
[a, b]_G$ is either empty or it is a singleton; in the latter
case, we denote the new element by $m(a, b)$.

 For a grid element $x$, the $C$-line through~$x$ is defined as
 \[
   \setm{a \in G}{a_D = x_D};
 \]
 symmetrically, we define the $D$-line through~$x$.

Note the following immediate consequences of (C1)--(C4):

 If $x$, $y \in L$ and $x \parallel y$, then
 \begin{enumerate}
 \item[(C5)] $x \mm y$, $x \jj y \in G$.
 \item[(C6)] $x \mm y = \ul x \mm \ul y$ and $x \jj y = \ol x
\jj \ol y$.
 \item[(C7)] $\ul x$ is on the $C$-line or on the $D$-line
through $x \mm y$; and symmetrically and dually.
 \end{enumerate}

The goal of this section is to prove the following result:

 \begin{theorem}\label{T:congrep}
 Let $L$ be a finite lattice satisfying conditions
\tup{(C1)--(C4)}. Then $L$ is $3$-modular.
 \end{theorem}

 \begin{proof}
 In this proof, we shall use the notation
 \[
   \vv<x, y, z>^{(n)} = \vv<x^{(n)}, y^{(n)}, z^{(n)}>,
 \]
 for a triple $\vv<x, y, z>$ in
$L$ and $n > 0$, and with this notation we can restate the
theorem:

\medskip

\emph{Let $L$ be a finite lattice satisfying conditions
\tup{(C1)--(C4)}. Then, for any triple $\vv<x, y, z>$ in $L$, the triple
$\vv<x^{(3)}, y^{(3)}, z^{(3)}>$ is balanced.}

\medskip

If $\vv<x, y, z>$ in $L$ is not an antichain, then $\vv<x^{(2)}, y^{(2)},
z^{(2)}>$ is balanced by Lemma~\ref{L:N5}. So from now on, we assume
that
 \begin{equation}
 \tag{A1} \text{$\vv<x, y, z>$ is an antichain.}
 \end{equation}

If $x \incomp y \mm z$, $y \incomp x \mm z$, $z \incomp x \mm
y$, then by (C5), $x^{(1)}$, $y^{(1)}$, $z^{(1)} \in G$, hence,
by Lemma~\ref{L:modular}, the triple $\vv<x^{(2)}, y^{(2)},
z^{(2)}>$ is balanced.  So by symmetry, we can assume:
 \begin{equation}
 \tag{A2}  y \mm z \leq x.
 \end{equation}
 Equivalently, $x^{(1)}=x$.

 If $x \mm z \leq y$ and $x \mm y \leq z$, then $\vv<x, y, z>$
is balanced. So we have two cases to consider:  $x \mm z \leq
y$, $ x \mm y \incomp z$ (see Figure~3) and $x \mm z \incomp
y$, $x \mm y \incomp z$ (see Figure~6).

\begin{figure}[htb]
\centerline{\includegraphics{Fig3.ill}}
\end{figure}

\tit{Case 1.} $x \mm z \leq y$ and $ x \mm y \incomp z$; see
Figure~3.

Note that the assumptions for Case 1 are symmetric in $x$ and
$y$.

Let $u = y \mm z$.  By (A2), $y \mm z \leq x \mm z$ and by the assumption
for Case 1, $x \mm z \leq y \mm z$, so
 \begin{equation}\label{E:xzyz}
   u = x \mm z = y \mm z \in G.
 \end{equation}

 Hence, $u \leq \ul z$.  We distinguish two subcases: $u = \ul
z$ and $u < \ul z$.

\tit{Case 1a.} $u = \ul z$.
 Obviously, $z \nin G$, because $z \in G$ would imply that $z =
\ul z = u \leq x$, contradicting the assumption (A1). So either
$z = m(u, \ol z)$ or $z = n(u, \ol z)$.  In~either case, $\ol z
\nleq x \mm y$, see Figures 4.1 and 4.2.  Let $L_C$ (resp,
$L_C'$) be the $C$-line through~$u$ (resp., through~$\ol z$).
By~symmetry and (C7), we can assume that $x \mm y$ is on~$L_C$.

Since $\ul x \mm \ul y = x \mm y$, it follows that either $\ul
x$ or $\ul y$ is on the line $L_C$. Since the assumptions are
symmetric in $x$ and $y$, we can assume that $\ul x$ is on the
line $L_C$; but $L_C'$ ``covers'' $L_C$ and $\ul y$ cannot
contain an element on $L_C'$ (that would contradict (A1)),
therefore, $\ul y$ also is on the line $L_C$.  We conclude that
$\ul x$ and $\ul y$ are comparable, so we can assume that $\ul
x \leq \ul y$ (and so $\ul x = x \mm y$).  In this case, $x =
\ul x$ leads to a contradiction with (A1), therefore, $\ul
x < x$. Also, $\ul x < \ul y$ because $\ul x = \ul y$ would
contradict either (A1) or~(C3). If $x = n(\ul x, \ol x)$ with
$\ul x$, $\ol x$ on $L_C$, then we cannot find room for $y$ by
(A1).

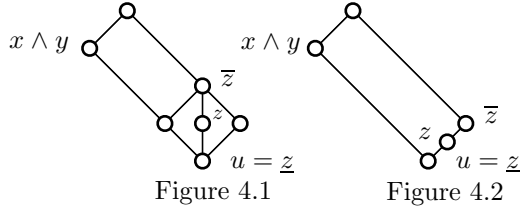
\begin{figure}[bth]
\centerline{\includegraphics{Fig4.ill}}
\end{figure}

So there are two possibilities for $x$:
 \begin{enumerate}
 \item $x = n(\ul x, \ol x)$ with $\ul x$ on $L_C$ and $\ol x$
on $L_C'$;
 \item $x = m(\ul x, \ol x)$ with $\ul x$ on $L_C$ and $\ol x$
on $L_C'$.
 \end{enumerate}

If (i) holds, then
\[
   \vv<x^{(1)}, y^{(1)}, z^{(1)}> = \vv<x, y,
\ol x > \text{\q and\q }x < \ol x,
 \]
 so $\vv<x^{(3)},
y^{(3)}, z^{(3)}>$ is balanced by Lemma~\ref{L:N5}.  If (ii)
holds, then
 \[
   \vv<x^{(1)}, y^{(1)}, z^{(1)}> = \vv<x, y, (x
\mm y) \jj z>
\]
 is balanced.

\begin{figure}[hbt]
\centerline{\includegraphics{Fig5.ill}}
\end{figure}

\tit{Case 1b.} $u < \ul z$.  Then $x^{(1)} = x$, $y^{(1)} = y$,
$z^{(1)} = (x \mm y) \jj z$, see Figure~5.  We~cannot have $\ul
x = \ul y$ ($ = x \mm y$);  indeed, then $x = \ul x$ or $y =
\ul y$ would contradict (A1); $\ul x < x$ and $\ul y < y$ would
contradict (C3).  It follows that by symmetry we can
assume without loss of generality that either $\ul x < \ul y$
or  $x \mm y < \ul x$, $x \mm y < \ul y$.  We consider these
cases separately.

\emph{Subcase I. $\ul x < \ul y$.} Obviously, $x \nin G$ because
$x \in G$ would imply that $x < y$, contradicting (A1). We
cannot have both $\ul x_C < \ul y_C$ and $\ul x_D < \ul y_D$
because this would again imply that $x < y$, contradicting (A1).

$\ul y$ and $z^{(1)}$ cannot be on the same line through $\ul
x$.  Indeed, if $\ul y$ and $z^{(1)}$ are, say, on the $D$-line
through $\ul x$, then $\ul x < \ul y < z^{(1)}$ (since $z^{(1)}
\leq \ul y$ would imply that $z < y$, contradicting (A1)).  So
$z^{(1)} = \ul x \jj \ol z = \ul y \jj \ol z$, which implies
that $\ul x \mm \ol z < \ul y \mm \ol z$ (since $G$ is
distributive). So $u < \ul y \mm \ul z$, contradicting $u = y
\mm z$.

Therefore, by symmetry, we can assume that $\ul x$ and $\ul y$
are on the $C$-line through $\ul x = x \mm y$ and $z^{(1)}$ is
on the $D$-line through $\ul x$.  Then $\ol x$ is not on the
$C$-line through $\ul x$ (this would contradict (A1)), so
either $x = n(\ul x, \ol x)$ with $\ul x_C = \ol x_C$, in which
case $\vv<x^{(3)}, y^{(3)}, z^{(3)}>$ is balanced (since $x =
x^{(1)} \leq z^{(1)}$) or $x = m(\ul x, \ol x)$, in which case
$\vv<x^{(1)}, y^{(1)}, z^{(1)}> = \vv<x, y, z^{(1)}>$ is
balanced.

\emph{Subcase II. $x \mm y < \ul x$ and $x \mm y < \ul y$.}
Let $L_C$ be the $C$-line through $x \mm y$ and let
$L_D$ be the $D$-line through $x \mm y$.  Since $z^{(1)} = (x
\mm y) \jj \ol z$, it follows that $z^{(1)}$ is on $L_C$ or on
$L_D$, say, on $L_D$.  Then by (A1), $\ul y < z^{(1)}$, which
contradicts \eqref{E:xzyz} as above, so this subcase cannot
occur.

This completes Case 1b and, therefore, Case 1.

\begin{figure}[hbt]
\centerline{\includegraphics{Fig6.ill}}
\end{figure}

\tit{Case 2.} $x \mm z \incomp y$ and $x \mm y \incomp z$.
This case is illustrated in Figure~6; the grey filled elements
are in $G$.  Define $u = (x \mm y) \jj (x \mm z)$.  The
elements $y \mm z$, $x \mm y$, $x \mm z$, and $u$ are distinct;
indeed, any equality would contradict with $x \mm z \incomp y$
or with $x \mm y \incomp z$. So $u$ is at least a ``prime
square'' above $y \mm z$.

We start with two observations:
 \begin{equation}\label{E:obs}
   \ul x = u\text{\q and\q }y^{(1)} \mm z^{(1)} = u.
 \end{equation}

Without loss of generality, we can assume that $x \mm y$ is on
the $C$-line through $y \mm z$ and that $x \mm z$ is on
the $D$-line through $y \mm z$.

If $u < \ul x$, then $(x \mm y)_C < \ul x_C$ or $(x \mm z)_D <
\ul x_D$, say, $(x \mm y)_C < \ul x_C$. It follows
that $\ul y_C < \ul x_C$ and so $\ol y_C \leq \ul x_C$.
Therefore, $y \leq x$, a~contradiction.

$y \mm z = \ul y \mm \ul z \leq u$ implies that $\ol y \mm \ol z
\leq u$ since $\ol y \mm \ol z$ is at most one ``prime square''
above $\ul y \mm \ul z$. Thus $(\ol y \jj u) \mm(\ol z
\jj u) = u$, as claimed.

Now this case is easy.  If $x = u$ or $x = m(\ul x, \ol x)$,
then $\vv<x^{(1)}, y^{(1)}, z^{(1)}>$ is balanced since $x =
x^{(1)}$.  If $x = n(\ul x, \ol x)$, then by symmetry we can
assume that $\ol x \leq z^{(1)}$ (since $u = \ul x$, $\ul x
\prec \ol x$, and $u < z^{(1)}$) and then $x^{(2)} = x$,
$y^{(2)} = y \jj \ol x$, $z^{(2)} = z^{(1)}$; hence
$\vv<x^{(2)}, y^{(2)}, z^{(2)}>$ is balanced, which completes
the proof.
 \end{proof}

Note that in Theorem~\ref{T:congrep}, ``$3$-modular'' cannot be changed
to ``$2$-modular''.  Indeed, let $A_0=\set{0, 1, 2}$ with $0 < 1 < 2$
and $A_1 = \set{0, 1}$  with $0 < 1$.  We take $a = \vv<1, 0>$, $b =
\vv<1, 1>$, and $L = G \uu \set{n(a, b)}$. Then $L$ satisfies
(C1)--(C4). Set $x = n(a, b)$, $y = \vv<2, 0>$, and $z = \vv<0, 1>$.
Then $x = x^{(1)} = x^{(2)} < x^{(3)} = \vv<1, 1>$.  So $L$ is not
$2$-modular.

\begin{corollary}\label{C:congrep}
 Every finite distributive lattice $D$ can be represented as the
congruence lattice of a finite planar $3$-modular lattice $L$.
 \end{corollary}

 \begin{proof}
 This is immediate by combining the representation theorem in G.
Gr\"atzer, H. Lakser, and E. T. Schmidt \cite{AGLS95} with
Theorem~\ref{T:congrep}.
 \end{proof}

\section{Discussion}\label{S:Discussion}
 \begin{problem}
 Which lattice varieties are closed under tensor product?
 \end{problem}

There are several ways to define what it means for a lattice
variety $\mathbf{V}$ to be ``closed under tensor product'':

\begin{enumerate}
\item If $A$ and $B$ are finite lattices, $A$, $B \in \mathbf{V}$,
then $A \otimes B
\in \mathbf{V}$.
\item If $A$ and $B \in \mathbf{V}$ are lattices with zero, $A$, $B
\in \mathbf{V}$, and
$A \otimes B$ is capped, then $A \otimes B \in \mathbf{V}$.
\item If $A$ and $B \in \mathbf{V}$ are lattices with zero, $A$, $B
\in \mathbf{V}$, and
$A \otimes B$ is a lattice, then $A \otimes B \in \mathbf{V}$.
\end{enumerate}

The trivial variety, $\mathbf{T}$, the variety of all
distributive lattices,
$\mathbf{D}$, and the variety of all lattices, $\mathbf{L}$, are
closed under tensor product under any of the three
interpretations. The problem, whether there are any more, is
open in all of its three variants.

It is easy to verify that $\mathbf{T}$ and $\mathbf{D}$ are the
only two finitely generated varieties that are closed under tensor
product (under any one of the interpretations).

 \begin{problem}
 Compute the modularity rank of $A[B]$, for small lattices $A$
and $B$?
 \end{problem}

The following examples were computed by B. Wolk:

\begin{example}\label{E:M3M4} The modularity rank of
$M_3[M_k]$ is $3$, for all
$k > 2$.
 \end{example}

Let $M_4 = \set{0,a,b,c,d,1}$. Then the computations can be
arranged in an array, as follows:

\smallskip

\begin{center}
\begin{tabular}{|| c | c | c | c ||}\hline $n$ & $p_n$ & $q_n$
& $r_n$ \\ \hline
$0$ & $\vv<b,c,a>$ & $\vv<b,a,d>$ & $\vv<a,0,c>$ \\ \hline $1$
& $\vv<b,c,a>$ &
$\vv<b,a,d>$ & $\vv<1,c,c>$ \\ \hline $2$ & $\vv<b,c,a>$ &
$\vv<1,1,1>$ &
$\vv<1,c,c>$ \\ \hline $3$ & $\vv<1,1,1>$ & $\vv<1,1,1>$ &
$\vv<1,1,1>$ \\
\hline
\end{tabular}
\end{center}

\smallskip

Among the 89,217 three-element antichains in $M_3[M_4]$, $936$
do not satisfy $\gm_2$ but they all satisfy $\gm_3$. More
generally, this is true in all the lattices
$M_3[M_k]$, $k > 2$, that is, they are all exactly $3$-modular.

\begin{example}
 Let $\mathbb{F}_7$ denote the lattice of
subspaces of the Fano plane. Then $M_3[\mathbb{F}_7]$ is not
$3$-modular.
\end{example}

Notation: the points are $1$, $2$, $3$, $4$, $5$, $6$, $7$;
the lines are $124$,
$235$, $346$,
$457$, $561$, $672$, $713$, and the plane is $PL$. The
following example shows that
$M_3[\mathbb{F}_7]$ is not $3$-modular.

\smallskip

\begin{center}
\begin{tabular}{|| c | c | c | c ||}\hline $n$ & $p_n$ & $q_n$
& $r_n$ \\ \hline
$0$ & $\vv<3,6,4>$ & $\vv<3,457,2>$ & $\vv<7,2,561>$ \\ \hline
$1$ &
$\vv<3,6,4>$ & $\vv<3,457,2>$ & $\vv<713,124,561>$ \\ \hline
$2$ &
$\vv<346,346,346>$ & $\vv<3,457,2>$ & $\vv<713,124,561>$ \\
\hline $3$ &
$\vv<346,346,346>$ & $\vv<713,457,672>$ & $\vv<713,124,561>$
\\ \hline $4$ &
$\vv<PL,346,346>$ & $\vv<713,457,672>$ & $\vv<713,124,561>$ \\
\hline
\end{tabular}
\end{center}

\smallskip

$M_3[\mathbb{F}_7]$ is too large for a complete search of the
type B. Wolk was conducting; it~has $1,090$ elements, and
around $190$ million three-element antichains.

Corollary~\ref{C:nonmod} show that the modularity rank of
$M_3[L]$ may be $\infty$ even if $L$ is modular.

\begin{problem} Let $K$ be a field, and let $V$ be a
$d$-dimensional vector space over $K$. What is the modularity
rank of $M_3[\LL(V)]$.
\end{problem}

\begin{problem}
 Is it possible to represent every finite distributive lattice
as the congruence lattice of a finite (planar) $2$-modular
lattice?
 \end{problem}

For $n > 1$, define $s(n)$ as the smallest integer so that
there is an exactly $n$-modular lattice of size $s(n)$.
Obviously, $s(2) = 5$, as realized by $N_5$.  The lattice
presented after the proof of Theorem~\ref{T:congrep} shows that
$s(3) = 7$.

\begin{problem}
 Determine the function $s(n)$.  Compute $s(n)$ for small
values of $n$.
 \end{problem}

 \begin{problem}
 Describe the free lattice with three generators over $\M_2$.
Is it finite?  What about $\tup{F}_{\M_n}(3)$?
 \end{problem}


\begin{thebibliography}{99}

\bibitem{DHW} A. Day, C. Herrmann, and R. Wille,
\emph{On modular lattices with four generators}, Algebra
Universalis \tbf{2} (1972), 317--323.

\bibitem{GLT} G. Gr\"atzer,
\emph{General Lattice Theory}, Pure and Applied Mathematics
\tbf{75}, Academic Press, Inc.\ (Harcourt Brace Jovanovich,
Publishers), New York-London; Lehrb\"ucher und Monographien
aus dem Gebiete der Exakten Wissenschaften, Mathematische
Reihe, Band 52. Birkh\"auser Verlag, Basel-Stuttgart; Akademie
Verlag, Berlin, 1978. xiii+381 pp.

\bibitem{AGLS95} G. Gr\"atzer, H. Lakser, and E. T. Schmidt,
\emph{Congruence lattices of small planar lattices}, Proc.\
Amer.\ Math.\ Soc.\ \tbf{123} (1995), 2619--2623.

\bibitem{GS95} G. Gr\"atzer and E. T. Schmidt,
\emph{A lattice construction and congruence-preserving
extensions}, Acta Math.\ Hungar.\ \tbf{66} (1995), 275--288.

\bibitem{GSa}
\bysame,
\emph{On the Independence Theorem of related structures for
modular (arguesian) lattices}, manuscript. Submitted for
publication in Studia Sci.\ Math.\ Hungar., March 1997.

\bibitem{GWa}
G. Gr\"atzer and F. Wehrung,
\emph{Proper congruence-preserving extensions of lattices},
AMS Abstract 97T-06-189. Acta Math.\ Hungar., to appear.

\bibitem{GWb}
\bysame,
\emph{Tensor products of semilattices with zero, revisited},
J. Pure Appl.\ Algebra, to appear.

\bibitem{GWc}
\bysame,
\emph{Tensor products and transferability of semilattices}, AMS
Abstract 97T-06-190.

\bibitem{MW73} A.~Mitchke and R. Wille,
\emph{Freie modulare Verb\"ande $FM(_DM_3)$.}  Proceedings of
the University of Houston Lattice Theory Conference (Houston,
Tex., 1973), pp. 383--396. Dept. Math., Univ. Houston,
Houston, Tex., 1973.

\bibitem{rQ85} R. W. Quackenbush,
\emph{Nonmodular varieties of semimodular lattices with a
spanning $M_3$.} Special volume on ordered sets and their
applications (L'Arbresle, 1982). Discrete Math.\ \tbf{53}
(1985), 193--205.

\bibitem{tS68} E. T. Schmidt,
\emph{Zur Charakterisierung der Kongruenzverb\"ande der
Verb\"ande}, Mat.\ \v Casopis Sloven.\ Akad.\ Vied.\ \tbf{18}
(1968), 3--20.

\bibitem{tS74}
\bysame,
\emph{Every finite distributive lattice is the congruence
lattice of a modular lattice}, Algebra Universalis \tbf{4}
(1974), 49--57.

\end{thebibliography}
\end{document}